\newtheorem{theorem}{Theorem}
\newtheorem{lemma}[theorem]{Lemma}
\newtheorem{corollary}[theorem]{Corollary}
\newtheorem{definition}[theorem]{Definition}

\newcommand{\qed}{{\hfill$\Box$}\\}

\newcommand{\iin}{\!\in\!}
\newcommand{\abs}[1]{\lvert#1\rvert}
\newcommand{\real}{\mathbb{R}}
\newcommand{\psm}{\Delta}
\newcommand{\norm}[1]{\lVert#1\rVert}
\newcommand{\normd}[1]{\lVert#1\rVert^\ast}
\newcommand{\normBL}[1]{\norm{#1}_{\mathsf{BL}}^\ast}

\newcommand{\msmeas}{\mu}
\newcommand{\nsmeas}{\nu}

\newcommand{\Xsp}{S}
\newcommand{\Ysp}{S^\prime}
\newcommand{\tMeas}{\mathsf{M_t}}
\newcommand{\tProb}{\mathsf{P_t}}

\newcommand{\Ub}{\mathsf{U_b}}
\newcommand{\BLip}{\mathsf{BLip_b}}
\newcommand{\Cb}{\mathsf{C_b}}
\newcommand{\Bo}{\mathsf{Bo}}
\newcommand{\ellone}{\ell_1}
\newcommand{\expect}{\mathsf{E}}

\newcommand{\mmap}{\varphi}

\newcommand{\uset}{\mathcal{U}}

\newcommand{\pmass}{\partial}

\documentclass{article}

\usepackage{amssymb,amsmath}

\title{Weak uniform structures on probability distributions}

\author{Jan Pachl  \\
Fields Institute\thanks{Paper written while the author was a visitor
at the Fields Institute.} \\
Toronto, Ontario, Canada}
\date{November 21, 2010 (version~2)}

\begin{document}
\maketitle

\begin{abstract}
In dealing with asymptotic approximation of possibly divergent
nets of probability distributions,
we are led to study uniform structures on the set of distributions.
This paper identifies a class of such uniform structures that may be considered to be
reasonable generalizations of the weak topology.
It is shown that all structures in the class yield the same notion of
asymptotic approximation for sequences (but not for general nets) of probability distributions.
\end{abstract}


\section{Introduction}
    \label{s:intro}

Traditionally limit theorems in probability theory are formulated in terms of
the \emph{weak topology} on probability distributions,
also known as the topology of convergence in distribution.
When a sequence of probability distributions converges in the weak topology,
its limit serves as a \emph{weak approximation} of the distributions in the sequence.

In more general limit theorems, a sequence of probability distributions,
not necessarily convergent, is asymptotically approximated by another sequence.
That leads to a natural question:
Are there essentially different ``reasonable'' notions of asymptotic approximation
that all reduce to the weak approximation in the special case of convergent sequences?

To deal with approximation for possibly divergent sequences we need more than a mere topology.
Asymptotic approximation is conveniently formulated in terms of uniform structures
(uniformities) on the space of probability distributions.
The question above then becomes:
What ``reasonable'' uniform structures on the space of probability distributions
are compatible with the weak topology?

Uniform structures on probability distributions on a metric space
were investigated by Dudley~\cite{Dudley1968dpm}\cite[11.7]{Dudley2002rap} and
D'Aristotile, Diaconis and Freedman~\cite{DAristotile1988omp}.
Davydov and Rotar~\cite{Davydov2009apd} showed that the previously neglected
uniformity defined by bounded uniformly continuous functions
is a reasonable uniformity in the above mentioned sense.

Uniform structures on probability distributions are also implicit in the rich literature
on probability metrics, surveyed by Dudley~\cite{Dudley1976pam},
Gibbs and Su~\cite{Gibbs2002cbp},
Rachev~\cite{Rachev1991pms} and Zolotarev~\cite{Zolotarev1997mts};
however, so far the focus in that area has been on quantitative results for convergent sequences,
not on general uniform structures or even the uniform structures
defined by probability metrics.

In this paper I offer further evidence that the uniformity proposed by Davydov and Rotar
is an appropriate analog of the weak topology,
at least when we deal with sequences (rather than general nets) of probability distributions.
To formalize the notion of ``reasonable'' uniform structures,
I formulate two natural properties that they should satisfy.
Although there is usually more than one uniformity satisfying those properties,
and therefore more than one notion of a ``weak'' asymptotic approximation for \emph{nets}
of probability distributions, they all agree on \emph{sequences}.

Version~1 of the paper was dated July 20, 2010.
Version~2 incorporates a suggestion from Ramon van Handel,
which considerably simplifies and improves the results in section~\ref{s:ustructures}.


\section{Preliminaries}
    \label{s:prelim}

Uniform structures (uniformities) and uniform spaces may be defined in several equivalent ways.
For the purposes of this paper, the most suitable definition is the one based
on pseudometrics~\cite[Ch.15]{Gillman1976rcf},
although definitions in terms of entourages~\cite{Bourbaki1940tg}
or uniform covers~\cite{Isbell1964us} would do as well.
Let $\uset$ be a uniform structure on a set $A$.
A net $\{a_\gamma\}_\gamma$ in $A$ is an \emph{asymptotic $\uset$-approximation}
of a net $\{b_\gamma\}_\gamma$ in $A$ iff
$\lim_\gamma \psm(a_\gamma , b_\gamma) = 0$
for every $\uset$-uniformly continuous pseudometric $\psm$.
Here a \emph{net} is a family indexed by a directed partially ordered set.
A \emph{sequence} is a family indexed by the totally ordered set $\{0,1,2,\dotsc\}$.

The set of real numbers is denoted $\real$.
When $\Xsp$ a metric space,
$\Cb(\Xsp)$ is the space of bounded real-valued continuous functions on $\Xsp$,
and $\Ub(\Xsp)$ is the space of bounded real-valued uniformly continuous functions on $\Xsp$.
The sup norm of a function $f\iin\Cb(\Xsp)$ is $\norm{f}_\Xsp := \sup_{x\in \Xsp} \abs{f(x)}$.

The \emph{Borel $\sigma$-algebra $\Bo(\Xsp)$} is the smallest $\sigma$-algebra
of subsets of $\Xsp$ containing all open sets in $\Xsp$.
In this paper, a \emph{measure on $\Xsp$} means a bounded signed measure on $\Bo(\Xsp)$.
A measure $\msmeas$ on $\Xsp$ is \emph{tight} if
\[
\abs{\msmeas}(A) = \sup \{ \abs{\msmeas}(C) \mid C \;\text{ is compact and }\; C \subseteq A \}
\]
for every $A\iin\Bo(\Xsp)$.
Tight measures are also known as Radon or regular measures.
The space of tight measures on $\Xsp$ is denoted $\tMeas(\Xsp)$, and
\[
\tProb(\Xsp) := \{ \msmeas\iin\tMeas(\Xsp) \mid \msmeas\geq 0 \;\text{ and }\; \msmeas(\Xsp)=1 \}
\]
is the space of tight probability measures on $\Xsp$.
The \emph{point mass at $x\iin \Xsp$} is denoted $\pmass_\Xsp(x)$;
that defines a mapping $\pmass_\Xsp \colon \Xsp \to \tProb(\Xsp)$.

The vector space duality $\langle\tMeas(\Xsp),\Cb(\Xsp)\rangle$ is defined by integration:
\[
\langle \msmeas, f \rangle := \int f \,\mathsf{d}\msmeas
\;\text{ for }\;\msmeas\iin\tMeas(\Xsp), f\iin\Cb(\Xsp).
\]
Since every $\msmeas\iin\tMeas(\Xsp)$ is uniquely determined by
the values $\langle \msmeas, f \rangle$, $f\iin\Ub(\Xsp)$,
integration defines also a vector space duality $\langle\tMeas(\Xsp),\Ub(\Xsp)\rangle$.
When $X$ is a real-valued random variable, $\expect(X)$ denotes the expected value of $X$.
If $f\iin\Cb(\Xsp)$ and $X$ is an $\Xsp$-valued random variable
with distribution $\msmeas\iin\tProb(\Xsp)$
then $\langle \msmeas, f \rangle = \expect(f(X))$.

When $\Xsp$ and $\Ysp$ are metric spaces and $\mmap\colon \Xsp\to \Ysp$ is a continuous mapping,
the linear mapping $\tMeas(\mmap)\colon\tMeas(\Xsp)\to\tMeas(\Ysp)$ is defined by
$\langle\tMeas(\mmap)(\msmeas),g\rangle:=\langle\msmeas,g\circ\mmap\rangle$
for $\msmeas\iin\tMeas(\Xsp)$ and $g\iin\Cb(\Ysp)$.
Clearly $\tMeas(\mmap)$ maps $\tProb(\Xsp)$ into $\tProb(\Ysp)$.

When $d$ is the metric of $\Xsp$ and $\msmeas\iin\tMeas(\Xsp)$, write
\begin{align*}
\BLip(d) & := \{ f\colon \Xsp \to \real \mid \norm{f}_\Xsp \leq 1 \;\text{ and }\;
\abs{f(x)-f(y)}\leq d(x,y) \;\text{ for all }\; x,y\iin \Xsp \} \\
\normd{\msmeas} & := \sup \left\{\langle \msmeas, f \rangle \mid f\iin\BLip(d) \right\} .
\end{align*}
By the following well-known lemma,
$\normd{\cdot}$ is a norm on the space $\tMeas(\Xsp)$.
This norm is equivalent to the norm $\normBL{\cdot}$ studied by
Dudley~\cite{Dudley1966cbm}\cite{Dudley1968dpm}.

\begin{lemma}
    \label{lem:approx}
Let $\Xsp$ be a metric space with metric $d$.
The space $\bigcup_{n=1}^\infty n\BLip(d)$ is $\norm{\cdot}_\Xsp$-dense in the space $\Ub(\Xsp)$.
\end{lemma}

\textbf{Proof}.
Take any $f\iin\Ub(\Xsp)$, $\varepsilon>0$.
There is $\theta>0$ such that if $x,y\iin \Xsp$ and $d(x,y)<\theta$
then $\abs{f(x)-f(y)}<\varepsilon$.
Choose an integer $n\geq \max(\norm{f}_\Xsp+\varepsilon,2\norm{f}_\Xsp/\theta)$
and define
\[
g(y) := \sup_{x \in \Xsp} \;
( \, f(x) - \varepsilon -  n \,d(x,y)\, ) \;\text{ for }\; y\iin \Xsp.
\]
Then $g\iin n\BLip(d)$ and $ f - \varepsilon \leq g \leq f $.
\qed

In addition to the $\normd{\cdot}$ topology,
consider also two weak topologies on $\tMeas(\Xsp)$:
\begin{itemize}
\item
The $\Cb(\Xsp)$-weak topology from the duality $\langle\tMeas(\Xsp),\Cb(\Xsp)\rangle$.
That is, the topology of simple convergence on the elements of $\Cb(\Xsp)$.
\item
The $\Ub(\Xsp)$-weak topology from the duality $\langle\tMeas(\Xsp),\Ub(\Xsp)\rangle$.
That is, the topology of simple convergence on the elements of $\Ub(\Xsp)$.
\end{itemize}
The restriction of the $\Cb(\Xsp)$-weak topology to
$\tProb(\Xsp)$ is called simply the \emph{weak topology}
(or the topology of convergence in distribution) in probability theory.

\begin{lemma}
    \label{lem:posweakt}
Let $\Xsp$ be a metric space.
The $\Cb(\Xsp)$-weak topology, the $\Ub(\Xsp)$-weak topology
and the $\normd{\cdot}$ topology coincide on $\tProb(\Xsp)$.
\end{lemma}

\textbf{Proof}.
LeCam~\cite[Lem.5]{LeCam1957cds} proves that the two weak topologies coincide on $\tProb(\Xsp)$.
The equivalence with the $\normd{\cdot}$ topology is proved by LeCam~\cite{LeCam1970ncc},
and by Dudley~\cite[Th.18]{Dudley1966cbm}.
\qed

The proof of Corollary~\ref{cor:prsequences} in section~\ref{s:ustructures}
relies on the following theorem,
which generalizes the Schur property~\cite[5.19]{Fabian2001fai}
of convergent sequences in $\ellone$.

\begin{theorem}
    \label{th:csequences}
Let $\Xsp$ be a metric space.
The $\Ub(\Xsp)$-weak topology and the $\normd{\cdot}$ topology on $\tMeas(\Xsp)$ have the same
compact sets, and therefore the same convergent sequences.
\end{theorem}

Theorem~\ref{th:csequences} appeared in~\cite{Pachl1979mfu}.
A simpler proof was found by Cooper and Schachermayer~\cite{Cooper1981umc}.
Variants of the theorem were proved by
van~Handel~\cite[B.1]{vanHandel2009uoh} (for measures on $\real^n$)
and Davydov and Rotar~\cite[Th.~4]{Davydov2009apd}.

\begin{theorem}
    \label{th:Dudley}
Let $\Xsp$ be a metric space and let $\msmeas_j, \nsmeas_j \iin \tProb(\Xsp)$ for $j=0,1,\dotsc$.
We have
\[
\lim_j \;\normd{\msmeas_j - \nsmeas_j} = 0
\]
if and only if there exist $\Xsp$-valued
random variables $X_j$, $Y_j$ such that the distribution of $X_j$ is $\msmeas_j$,
the distribution of $Y_j$ is $\nsmeas_j$, and $\lim_j d(X_j,Y_j) = 0$ almost surely.
\end{theorem}

Theorem~\ref{th:Dudley} is proved by Dudley~\cite[11.7.1]{Dudley2002rap}.
(In~\cite{Dudley2002rap}, $\beta$ denotes the metric of the norm $\normBL{\cdot}$,
which is equivalent to $\normd{\cdot}$.)

When $E$ is a locally convex vector space,
the topology of $E$ is defined by the family of continuous seminorms~\cite[II.4]{Schaefer1971tvs}.
The \emph{additive uniformity} on any subset of $E$
is the uniformity induced by the pseudometrics of the form $(x,y)\mapsto s(x-y)$, $x,y\iin E$,
where $s$ is a continuous seminorm on $E$.
Each of the three topologies on $\tMeas(\Xsp)$ defined above makes $\tMeas(\Xsp)$
into a locally convex space.
The corresponding additive uniformities on subsets of $\tMeas(\Xsp)$
(in particular, on $\tProb(\Xsp)$)
will be referred to as the $\normd{\cdot}$ uniformity,
the $\Cb(\Xsp)$-weak uniformity and the $\Ub(\Xsp)$-weak uniformity.

\begin{lemma}
    \label{lem:twouniform}
Let $\Xsp$ be a metric space.
On $\tProb(\Xsp)$, the $\normd{\cdot}$ uniformity is finer than the $\Ub(\Xsp)$-weak uniformity.
\end{lemma}

\textbf{Proof}.
Follows from Lemma~\ref{lem:approx}.
\qed


\section{Uniform structures on $\tProb(\Xsp)$}
    \label{s:ustructures}

In this section I describe a class of uniform structures on $\tProb(\Xsp)$
that appear to be reasonable candidates for extending to divergent nets
the notion of weak approximation on convergent nets.
Of course, it is highly subjective and context-dependent what uniformities
should be considered ``reasonable'' for this purpose.
I replace that subjective notion by the two properties in Definition~\ref{def:props}.

\begin{lemma}
    \label{lem:compunif}
Let $I_n := [-n,n] \subseteq \real$ with the standard metric.
There is a unique uniformity on $\tProb(I_n)$ compatible with
the $\Cb(I_n)$-weak topology.
This unique uniformity coincides with the $\normd{\cdot}$ uniformity,
the $\Cb(I_n)$-weak uniformity and the $\Ub(I_n)$-weak uniformity on $\tProb(I_n)$.
\end{lemma}

\textbf{Proof}.
Since $I_n$ is compact, so is $\tProb(I_n)$ with the $\Cb(I_n)$-weak topology.
It follows that there is only one uniformity on $\tProb(I_n)$
compatible with the $\Cb(I_n)$-weak topology~\cite[II.24]{Isbell1964us},
and by Lemma~\ref{lem:posweakt} it coincides with the three uniformities
listed.
\qed

Property (A1) in the next definition and the resulting simplification of the proofs that follow
were suggested by R. van Handel.

\begin{definition}
    \label{def:props}
Let $\Xsp$ be a metric space with metric $d$.
Consider the following properties of a~uniform structure $\uset$ on $\tProb(\Xsp)$.
\begin{itemize}
\item[\textrm{(A1)}]
If $X_j$ and $Y_j$, $j=0,1,\dotsc$, are $\Xsp$-valued random variables
with distributions $\msmeas_j\iin\tProb(\Xsp)$ and $\nsmeas_j\iin\tProb(\Xsp)$
respectively and if $\lim_j d(X_j,Y_j) = 0$ almost surely
then the sequence $\{\msmeas_j\}_j$ is an asymptotic $\uset$-approximation
for the sequence $\{\nsmeas_j\}_j$.
\item[\textrm{(A2)}]
For every interval $I_n:=[-n,n]\subseteq\real$, $n=1,2,\dotsc$,
if a mapping $\mmap\colon \Xsp \to I_n$ is uniformly continuous then so is the mapping
$\tMeas(\mmap)$ from $\tProb(\Xsp)$ with $\uset$
to $\tProb(I_n)$ with the unique uniformity in Lemma~\ref{lem:compunif}.
\end{itemize}
\end{definition}

Property (A1) states that asymptotic approximation (understood almost surely) for sequences
of $\Xsp$-valued random variables implies asymptotic approximation
for the corresponding distributions.

(A2) is a functorial property of the assignment $X\mapsto(\tProb(\Xsp),\uset)$
for uniformly continuous mappings $\mmap\colon \Xsp \to I_n$,
assuming that on convergent nets in $\tProb(I_n)$
the asymptotic approximation agrees with the weak approximation.

As a special case of (A1), if $\{x_j\}_j$ and $\{y_j\}_j$
are two sequences of points in $\Xsp$ such that $\lim_j d(x_j,y_j) = 0$ then
$\{\pmass(x_j)\}_j$ is an asymptotic $\uset$-approximation for $\{\pmass(y_j)\}_j$.
Consider the points $x_j:=j$ and $y_j=j+1/j$, $j=1,2,\dotsc$, in $\real$.
There is a function $f\iin\Cb(\real)$ such that $f(x_j)-f(y_j)=1$ for all $j$;
it follows that the $\Cb(\real)$-weak uniformity on $\tProb(\real)$ does not satisfy (A1).

The other two uniformities from section~\ref{s:prelim} have properties (A1) and (A2):

\begin{theorem}
    \label{th:prproperties}
For every metric space $\Xsp$
the $\Ub(\Xsp)$-weak uniformity
and the $\normd{\cdot}$ uniformity on $\tProb(\Xsp)$ have properties (A1) and (A2).
\end{theorem}

\textbf{Proof}.
Let $\msmeas_j,\nsmeas_j\iin\tProb(X)$, $j=0,1,\dotsc$,
be the distributions of $\Xsp$-valued random variables $X_j$, $Y_j$
such that $\lim_j d(X_j,Y_j) = 0$ almost surely.
Then
\begin{gather*}
\abs{\langle\msmeas_j,f\rangle - \langle\nsmeas_j,f\rangle}
= \abs{\expect(f(X_j)) - \expect(f(Y_j))}
\leq \expect\left(\abs{f(X_j)-f(Y_j)}\right)            \\
\normd{\msmeas_j - \nsmeas_j}
= \sup_{f\in\BLip(d)} \abs{\langle\msmeas_j,f\rangle - \langle\nsmeas_j,f\rangle}
\leq \expect( 2 \wedge d(X_j,Y_j))
\end{gather*}
and if $\lim_j d(X_j,Y_j) = 0$ almost surely then  $\lim_j \;\normd{\msmeas_j - \nsmeas_j} = 0$.
Hence the $\normd{\cdot}$ uniformity has property (A1),
and so does the $\Ub(X)$-weak uniformity by Lemma~\ref{lem:twouniform}.

(A2) follows from the definition of $\tMeas(\mmap)$ and Lemma~\ref{lem:compunif}.
\qed

Next I prove that every uniformity $\uset$ satisfying (A1) and (A2)
is between the two uniformities in Theorem~\ref{th:prproperties}.

\begin{theorem}
    \label{th:prpropbetween}
Let $X$ be a metric space with metric $d$, and $\uset$ a uniform structure on $\tProb(X)$.
\begin{enumerate}
\item
If $\uset$ has property (A1) then
it is coarser than the $\normd{\cdot}$ uniformity on $\tProb(X)$.
\item
If $\uset$ has property (A2)
then it is finer than the $\Ub(X)$-weak uniformity on $\tProb(X)$.
\end{enumerate}
\end{theorem}

\textbf{Proof}.
To prove part~1, assume that $\uset$ is not coarser
than the $\normd{\cdot}$ uniformity on $\tProb(X)$.
That means that there are a $\uset$-uniformly continuous pseudometric $\psm$, $\varepsilon>0$
and $\msmeas_j,\nsmeas_j\iin\tProb(X)$ for $j=0,1,\dotsc$
such that $\lim_j \;\normd{\msmeas_j-\nsmeas_j} =0$ and
$\psm(\msmeas_j,\nsmeas_j)\geq\varepsilon$ for all $j$.
By Theorem~\ref{th:Dudley} there are $\Xsp$-valued random variables
$X_j$, $Y_j$
whose distributions are $\msmeas_j,\nsmeas_j$ respectively and such that
$ \lim_j \;d(X_j,Y_j) = 0 $ almost surely.
Thus $\uset$ does not have property (A1).

2. Assume that $\uset$ has property (A2) and take any $f\iin\Ub(X)$.
Choose $n\iin\{1,2,\dotsc\}$ for which $\norm{f}_X\leq n$,
so that $f$ maps $X$ into $I_n:=[-n,n]$.

Let $g\iin\Ub(I_n)$ be the function $g:x\mapsto x$.
By (A2) and Lemma~\ref{lem:compunif}
the mapping $\tMeas(f)$ is uniformly continuous from $\tProb(X)$
with $\uset$ to $\tProb(I_n)$ with the $\Ub(I_n)$-weak uniformity.
Thus the mapping
\[
\msmeas\mapsto \langle\tMeas(f)(\msmeas),g\rangle
= \langle \msmeas, g\circ f \rangle\ = \langle \msmeas, f \rangle
\]
from $\tProb(X)$ with $\uset$ to $\real$ is uniformly continuous.
That proves that $\uset$ is finer than the $\Ub(X)$-weak uniformity.
\qed

\begin{corollary}
    \label{cor:compatible}
Let $X$ be a metric space with metric $d$, and $\uset$ a uniform structure on $\tProb(X)$.
If $\uset$ has properties (A1) and (A2) then it is compatible
with the $\Cb(X)$-weak topology on $\tProb(X)$.
\end{corollary}

\textbf{Proof}.
Apply Theorem~\ref{th:prpropbetween} and Lemma~\ref{lem:compunif}.
\qed

\begin{corollary}
    \label{cor:prsequences}
Let $\Xsp$ be a metric space
and $\uset$ a uniform structure on $\tProb(\Xsp)$.
If $\uset$ has properties (A1) and (A2) then the following statements
are equivalent for any two sequences $\{\msmeas_j\}_j$ and $\{\nsmeas_j\}_j$
in $\tProb(\Xsp)$.
\begin{enumerate}
\item[\textit{(i)}]
$\lim_j \,(\langle\msmeas_j,f\rangle - \langle\nsmeas_j,f\rangle) = 0$ for every $f\iin\Ub(\Xsp)$.
\item[\textit{(ii)}]
The sequence $\{\msmeas_j\}_j$ is an asymptotic $\uset$-approximation
of the sequence $\{\nsmeas_j\}_j$\,.
\item[\textit{(iii)}]
$\lim_j \,\normd{\msmeas_j - \nsmeas_j}=0$.
\item[\textit{(iv)}]
There exist $\Xsp$-valued random variables $X_j$, $Y_j$, $j=0,1,\dotsc$,
whose distributions are $\msmeas_j$ and
$\nsmeas_j$ and such that $\lim_j \,d(X_j,Y_j)=0$ almost surely.
\end{enumerate}
\end{corollary}

\textbf{Proof}.
The implication (i)$\Rightarrow$(iii) follows from Theorem~\ref{th:csequences},
(iii)$\Rightarrow$(ii)$\Rightarrow$(i) from Theorem~\ref{th:prpropbetween},
and (iii)$\Leftrightarrow$(iv) from Theorem~\ref{th:Dudley}.
\qed

The equivalence of (i) and (iii) in the corollary was derived by
Davydov and Rotar~\cite{Davydov2009apd} from their variant of Theorem~\ref{th:csequences}.
By~\cite[11.7.1]{Dudley2002rap}, condition (iii) is also equivalent to each of
the following:
\begin{itemize}
\item
$\lim_j \,\rho(\msmeas_j,\nsmeas_j)=0$, where
$\rho$ is the L\'{e}vy--Prokhorov metric.
\item
There exist $\Xsp$-valued random variables $X_j$, $Y_j$, $j=0,1,\dotsc$,
whose distributions are $\msmeas_j$ and
$\nsmeas_j$ and such that $\lim_j \,d(X_j,Y_j)=0$ in probability.
\end{itemize}


\section{Concluding remarks}
    \label{s:concl}

By Corollary~\ref{cor:prsequences},
all uniformities on $\tProb(\Xsp)$ that satisfy (A1) and (A2) yield the same notion
of asymptotic approximation for \emph{sequences} of tight probability measures,
namely the approximation defined by the $\Ub(\Xsp)$-weak and $\normd{\cdot}$ uniformities.

When we look beyond sequences and deal with general \emph{nets} of probability distributions,
there are multiple notions of asymptotic approximation,
even for the uniformities with properties (A1) and (A2).
For example, the $\normd{\cdot}$ uniformity on $\tProb(\real)$ is strictly finer than
the $\Ub(\real)$-weak uniformity,
and thus there are two nets $\{\msmeas_\gamma\}_\gamma$ and $\{\nsmeas_\gamma\}_\gamma$
in $\tProb(\real)$ such that $\{\nsmeas_\gamma\}_\gamma$ is an asymptotic
$\Ub(\real)$-weak approximation of $\{\msmeas_\gamma\}_\gamma$
but not an asymptotic $\normd{\cdot}$ approximation.
\\

\textbf{Acknowledgement.}
I wish to thank David Fremlin, Ramon van Handel and B{\'a}lint Vir{\'a}g for their comments.


\end{document}